\documentclass[a4paper,reqno,11pt,oneside]{amsart}
\usepackage{amsmath,geometry}
\usepackage{graphicx}
\usepackage{mathrsfs}
\usepackage{amssymb,amsmath, amsfonts, amsthm}
\usepackage{latexsym}
\usepackage{color} 
\usepackage[dvips]{epsfig}
\usepackage[colorlinks]{hyperref}
\usepackage{bm}

\theoremstyle{plain}
\newtheorem{definition}{Definition}[section]
\newtheorem{thm}[definition]{Theorem}

\newtheorem{lem}[definition]{Lemma}

\newtheorem{rem}[definition]{Remark}

\begin{document}
\title[A note on the multiplicity of eigenvalues]{A note on the persistence of multiplicity of eigenvalues of fractional Laplacian
under perturbations}
 
\author{Marco Ghimenti}
\address[Marco Ghimenti]{Dipartimento di Matematica
Universit\`a di Pisa
Largo Bruno Pontecorvo 5, I - 56127 Pisa, Italy}
\email{marco.ghimenti@unipi.it }

\author{Anna Maria Micheletti}
\address[Anna Maria Micheletti]{Dipartimento di Matematica
Universit\`a di Pisa
Largo Bruno Pontecorvo 5, I - 56127 Pisa, Italy}
\email{a.micheletti@dma.unipi.it }

\author{Angela Pistoia}
\address[Angela Pistoia] {Dipartimento SBAI, Universit\`{a} di Roma ``La Sapienza", via Antonio Scarpa 16, 00161 Roma, Italy}
\email{angela.pistoia@uniroma1.it}

\thanks{The first author is partially supported by the MIUR Excellence Department Project awarded to the Department of Mathematics, University of Pisa, CUP I57G22000700001 and by the PRIN 2022 project 2022R537CS \emph{$NO^3$ - Nodal Optimization, NOnlinear elliptic equations, NOnlocal geometric problems, with a focus on regularity}, founded by the European Union - Next Generation EU. The second author is partially supported by the MUR-PRIN-20227HX33Z ``Pattern formation in nonlinear phenomena''. 
The research of the authors is partially supported by the GNAMPA project 2024: ``Problemi di doppia
curvatura su variet\`a a bordo e legami con le EDP di tipo ellittico''}

\begin{abstract}
We consider the eigenvalues problem for the the fractional Laplacian
$(-\Delta)^{s}$, $s\in(0,1)$, in a bounded domain $\Omega$ with
Dirichlet boundary condition. A recent result  (see \emph{Generic properties of eigenvalues of the fractional Laplacian} by Fall, Ghimenti,  Micheletti and Pistoia, 
CVPDE (2023)) states 
 that under generic small perturbations of the coefficient of
the equation or of the domain $\Omega$ all the eigenvalues are simple.
In this paper we give a condition for which a perturbation of the
coefficient or of the domain preserves the multiplicity of a given
eigenvalue. Also, in the case of an eigenvalue of multiplicity $\nu=2$
we prove that the set of perturbations of the coefficients which preserve
the multiplicity is a smooth manifold of codimension $2$ in $C^{1}(\Omega)$.
\end{abstract}


\maketitle
\section{Introduction}

In the last decade, there has been a great deal of interest in using the fractional Laplacian to model diverse physical phenomena.
We refer the readers to    Di Nezza, Palatucci  and Valdinoci’s survey paper \cite{DPV} for a detailed exposition of the function spaces involved in the analysis of the operator and  to the   recent Ros-Oton's expository paper \cite{RO} for a list of  results on Dirichlet problems on bounded domains. \\

 In this paper we will focus on the eigenvalue problem

 \begin{equation}
\left\{ \begin{array}{cc}
(-\Delta)^{s}\varphi_{s}=\lambda\varphi_{s} & \text{ in }\Omega\\
\\
\varphi_{s}=0 & \text{ in }\Omega^{c}=\mathbb{R}^{n}\smallsetminus\Omega
\end{array}\right.,\label{eq:Pb-2}
\end{equation}
where $(-\Delta)^{s}$ for $0<s<1$ denotes the fractional Laplacian
and $\Omega$ is a $C^{1,1}$ bounded domain in $\mathbb{R}^{n}$,
with $n>2s$. 

In a weak sense, problem \eqref{eq:Pb-2} can be formulated as follows. We consider the
space 
\[
\mathcal{H}_{0}^{s}(\Omega):=\left\{ u\in H^{s}(\mathbb{R}^{n})\ :\ u\equiv0\text{ on }\Omega^{c}\right\} ,
\]
where 
\[
H^{s}(\mathbb{R}^{n}):=\left\{ u\in L^{2}(\mathbb{R}^{n})\ :\ \frac{u(x)-u(y)}{|x-y|^{\frac{n}{2}+s}}\in L^{2}(\mathbb{R}^{n}\times\mathbb{R}^{n})\right\} ,
\]
and the quadratic form defined on  $\mathcal{H}_{0}^{s}(\Omega)$ by
\[
(u,v)\mapsto\mathcal{E}_{s}^{\Omega}(u,v)=\mathcal{E}(u,v):=\frac{C_{n,s}}{2}\int_{\mathbb{R}^{n}}\int_{\mathbb{R}^{n}}\frac{(u(x)-u(y))(v(x)-v(y))}{|x-y|^{n+2s}}dxdy
\]
Then, we say that $\varphi_{s}\in\mathcal{H}_{0}^{s}(\Omega)$ is
an eigenfunction corresponding to the eigenvalue $\lambda_{s}$ iff
\[
\mathcal{E}(\varphi_{s},v)=\lambda_{s}\int_{\mathbb{R}^{n}}\varphi_{s}vdx\ \ \forall v\in\mathcal{H}_{0}^{s}(\Omega).
\]
It is well known (see, for instance, \cite{BRS}) that (\ref{eq:Pb-2})
admits an ordered sequence of eigenvalues 
\[
0<\lambda_{1,s}<\lambda_{2,s}\le\lambda_{3,s}\le\dots\le\lambda_{1,s}\le\dots\rightarrow+\infty.
\]

Since the first eigenvalue is strictly positive, we can also endow
$\mathcal{H}_{0}^{s}(\Omega)$ with the norm
\[
\|u\|_{\mathcal{H}_{0}^{s}(\Omega)}^{2}=\|u\|_{L^{2}(\Omega)}^{2}+\mathcal{E}(u,u).
\]

We refer to \cite{frank} and the references therein for a review of  results on eigenvalues of fractional Laplacians and fractional Schrödinger operators.
\\

In the recent paper \cite{FGMP}, Fall, Ghimenti, Micheletti and Pistoia     prove that there exist arbitrarily
small perturbations of the domain or arbitrarily small perturbations
of the coefficient of the linear terms for which all the eigenvalues of problems 
\begin{equation}
\left\{ \begin{array}{cc}
(-\Delta)^{s}\varphi_{s}+a(x)\varphi_{s}=\lambda\varphi_{s} & \text{ in }\Omega\\
\\
\varphi_{s}=0 & \text{ in }\Omega^{c}=\mathbb{R}^{n}\smallsetminus\Omega
\end{array}\right.\label{eq:Pb-coeff1}
\end{equation}
and 
\begin{equation}
\left\{ \begin{array}{cc}
(-\Delta)^{s}\varphi_{s}=\lambda a(x)\varphi_{s} & \text{ in }\Omega\\
\\
\varphi_{s}=0 & \text{ in }\Omega^{c}=\mathbb{R}^{n}\smallsetminus\Omega
\end{array}\right.\label{eq:Pb-coeff2}
\end{equation}
are simple.\\

In this paper we want to study the structure of the set of  perturbations of
the coefficients or  of the domain which \emph{preserve} the
multiplicity of the eigenvalues\footnote{This question was raised by the anonymous referee. We wish to thank
them for their interesting suggestion}. \\

Our first result deals with  the perturbation of the coefficients.
\begin{thm}
\label{thm:coeff}Let $\lambda_{0}$ be an eigenvalue for Problem
(\ref{eq:Pb-coeff1}) (respectively Problem (\ref{eq:Pb-coeff2}))
with multiplicity $\nu>1$, and let $\varphi_{1},\dots,\varphi_{\nu}$
be an $L^{2}$-orthonormal basis for the eigenspace relative to $\lambda_{0}$.
Assume that $a\in C^{1}(\Omega)$ and $\min_{\bar{\Omega}}a>0$ or
$\|a\|_{C^{1}(\Omega)}$ small (resp. assume that $a\in C^{1}(\Omega)$
and $\min_{\bar{\Omega}}a>0$). Let $b\in C^{1}(\Omega)$ be sufficiently
small and consider the functionals
\begin{equation}
b\mapsto\gamma_{ij}(b):=\int_{\Omega}b\varphi_{i}\varphi_{j},\ \ i,j=1,\dots,\nu\label{eq:spezz-coeff}
\end{equation}
 Then the set $\mathscr{I}$ of the $b$'s close to $0$ in $C^1(\Omega)$ such that the perturbed
problem
\begin{eqnarray*}
(-\Delta)^{s}\varphi+\left(a(x)+b(x)\right)\varphi=\lambda\varphi\text{ in }\Omega, &  & \varphi=0\text{ in }\Omega^{c}
\end{eqnarray*}
(respectively $(-\Delta)^{s}\varphi+\varphi=\lambda\left(a(x)+b(x)\right)\varphi\text{ in }\Omega,\varphi=0\text{ in }\Omega^{c}$)
admits an eigenvalue $\lambda_{b}$ close to ${\lambda_0}$ of the same multiplicity
$\nu$ is a subset of 
\[
\mathscr{H}:=\left\{ b\in C^{1}(\Omega)\ :\ \gamma_{ij}(b)=0\text{ for }i\neq j,\ \gamma_{11}(b)=\gamma_{22}(b)=\dots=\gamma_{\nu\nu}(b)\right\} .
\]
In addition, if the map 
\begin{align*}
G: & C^{1}(\Omega)\rightarrow L(\mathbb{R}^{\nu},\mathbb{R}^{\nu})\\
G(b) & =\left(\gamma_{ij}(b)\right)_{ij}
\end{align*}
is such that the span of $G(b)$ and the Identity map gives all the
$\nu\times\nu$ symmetric matrices, then the set $\mathscr{I}$ is
a manifold in $C^{1}(\Omega)$ of codimension $\frac{\nu(\nu+1)}{2}-1$. \\
In particular the last claim holds if $\lambda_{0}$ is an eigenvalue
of multiplicity $\nu=2$. 
\end{thm}

Our second result deals with the  perturbation of the domain.
\begin{thm}
\label{thm:dom}Let $\lambda_{0}$ be an eigenvalue for Problem (\ref{eq:Pb-2})
with multiplicity $\nu>1$, and let $\varphi_{1},\dots,\varphi_{\nu}$
be an $L^{2}$-orthonormal basis for the eigenspace relative to $\lambda_{0}$.
Let $\psi\in C^{1}(\mathbb{R}^{n},\mathbb{R}^{n})$ sufficiently small
and consider the functionals
\begin{equation}
\psi\mapsto\gamma_{ij}(\psi):=\int_{\partial\Omega}\frac{\varphi_{i}}{\delta^{s}}\frac{\varphi_{j}}{\delta^{s}}\psi\cdot N,\ \ i,j=1,\dots,\nu\label{eq:spezz-dom}
\end{equation}
 where $\delta(x)=\mathrm{dist}(x,\mathbb{R}^{n}\smallsetminus\Omega)$
and $N$ is the exterior normal of $\partial\Omega$

Then the set $\mathscr{I}$ of the $\psi$'s close to $0$ in $C^{1}(\mathbb{R}^{n},\mathbb{R}^{n})$ such that the problem
\begin{eqnarray*}
(-\Delta)^{s}\varphi+\varphi=\lambda\varphi\text{ in }\Omega_{\psi}, &  & \varphi=0\text{ in }\Omega_{\psi}^{c}
\end{eqnarray*}
in the perturbed domain $\Omega_{\psi}$ admits an eigenvalue $\lambda_{\psi}$ close to $ {\lambda_0}$
of the same multiplicity $\nu$ is a subset of 
\[
\mathscr{H}:=\left\{ \psi\in C^{1}(\mathbb{R}^{n},\mathbb{R}^{n})\ :\ \gamma_{ij}(\psi)=0\text{ for }i\neq j,\ \gamma_{11}(\psi)=\gamma_{22}(\psi)=\dots=\gamma_{\nu\nu}(\psi)\right\} .
\]
In addition, if the map 
\begin{align*}
G: & C^{1}(\mathbb{R}^{n},\mathbb{R}^{n})\rightarrow L(\mathbb{R}^{\nu},\mathbb{R}^{\nu})\\
G(\psi) & =\left(\gamma_{ij}(\psi)\right)_{ij}
\end{align*}
is such that the span of $G(\psi)$ and the Identity map gives all
the $\nu\times\nu$ symmetric matrices, then the set $\mathscr{I}$
is a manifold in $C^{1}(\mathbb{R}^{\nu},\mathbb{R}^{\nu})$ of codimension
$\frac{\nu(\nu+1)}{2}-1$. 
\end{thm}

The proof of our  results follows the strategy developed by Micheletti and Lupo in
\cite{ML,DL}, where an abstract transversality result is applied
to a second order elliptic operator under the effect of the perturbations
of the domain. The application of the abstract theorem in the case
of multiplicity $\nu=2$, which gives a concrete example of $\mathscr{I}$
being a manifold, relies on the unique continuation property. For
nonlocal problem this property has been proved only in particular
setting, and it is a challenging field of research. That is why
we can prove that $\mathscr{I}$ is a manifold only for problem (\ref{eq:Pb-coeff1})
and (\ref{eq:Pb-coeff2}). Actually in Remark \ref{dom} we will point out what it would be necessary to complete the proof in the case of
the perturbation of the domain  for problem \eqref{eq:Pb-2}.\\

The paper is organized as follows. Firstly we recall the abstract
transversality theorem. Then we prove the result for Problem (\ref{eq:Pb-coeff1})
and we sketch the proof for Problem (\ref{eq:Pb-coeff2}), concluding
the proof of Thm \ref{thm:coeff}. In the last section we prove Thm
\ref{thm:dom}.

\section{The abstract transversality result}

We recall here an abstract result which holds in a Hilbert space $X$
endowed with a scalar product $<\cdot,\cdot>_{X}$ for a selfadjoint
compact operator $T_{b}:X\rightarrow X$ depending smoothly on a parameter
$b$ which is defined in some Banach space $B$. If $T_{0}$ admits
an eigenvalue $\bar{\lambda}$ with multiplicity $\nu>1$, we provide
a characterization for the set $\mathscr{I}$ of parameter $b$ for
which $T_{b}$ has an eigenvalue $\lambda_{b}$ near $\bar{\lambda}$
which maintains the same multiplicity $\nu$ is a manifold in $B$.
In addition the result gives a sufficient condition which ensure that
$\mathscr{I}$ is a smooth sub-manifold of $B$.
\begin{thm}
\label{thm:astratto-1}Let $T_{b}:X\rightarrow X$ be a selfadjoint
compact operator which depends smoothly on a parameter $b$ belonging
to a real Banach space $B$. Let $T_{0}=T$ and $T_{b}$ be Frechet
differentiable in $b=0$. Let Let $x_{1}^{0},\dots,x_{\nu}^{0}$ be
an orthonormal base for the eigenspace relative to $\bar{\lambda}$.
If $b$ is sufficiently small to ensure $T_{b}\in\tilde{M}$, and
consider the functionals
\begin{equation}
b\mapsto\gamma_{ij}(b):=<T'(0)[b]x_{j}^{0},x_{i}^{0}>_{X},\ \ i,j=1,\dots,\nu\label{eq:spezzamentoastratto-1}
\end{equation}
 Then the set $\mathscr{I}$ of the $b$'s close to  $0$ in $X$ such that $T_{b}$
admits an eigenvalue $\lambda_{b}\sim\bar{\lambda}$ of the same multiplicity
$\nu$ is a subset of 
\[
\mathscr{H}:=\left\{ b\in B\ :\ \gamma_{ij}(b)=0\text{ for }i\neq j,\ \gamma_{11}(b)=\gamma_{22}(b)=\dots=\gamma_{\nu\nu}(b)\right\} .
\]
In addition, if the map 
\begin{align*}
G: & B\rightarrow L(\mathbb{R}^{\nu},\mathbb{R}^{\nu})\\
G(b) & =\left(\gamma_{ij}(b)\right)_{ij}
\end{align*}
is such that the span of $G(b)$ and the Identity map gives all the
$\nu\times\nu$ symmetric matrices, then the set $\mathscr{I}$ is
a manifold in $B$ of codimension $\frac{\nu(\nu+1)}{2}-1$. 
\end{thm}

The first part ot Theorem \ref{thm:astratto-1} was firstly proved
in \cite{Mi76}. A sketched version of the proof can be found also
in \cite{FGMP}, since condition (\ref{eq:spezzamentoastratto-1})
was the main tool to prove that the eigenvalues for fractional laplacian
are generically simple under perturbation of the domain or of the
coefficients. The proof of the second part can be found in \cite[,Th. 1]{DL}. 

\section{The case of Problem (\ref{eq:Pb-coeff1})}

We consider on n $\mathcal{H}_{0}^{s}(\Omega)$ the quadratic form
\[
\mathcal{B}^{a}(u,v)=\mathcal{E}(u,v)+\int_{\mathbb{R}^{n}}au^{2}dx.
\]
Since $\min_{\bar{\Omega}}a>0$ or $\|a\|_{C^{1}(\Omega)}$ small
, $\mathcal{B}^{a}(u,v)$ is a positive definite scalar product, and
we can consider on $\mathcal{H}_{0}^{s}(\Omega)$ the equivalent norm
\begin{equation}
\|u\|_{\mathcal{H}_{0}^{s}(\Omega)}^{2}=\mathcal{B}^{a}(u,u)=\mathcal{E}(u,u)+\int_{\mathbb{R}^{n}}au^{2}dx.\label{eq:scalar}
\end{equation}

Given the continuous and compact embedding $i:\mathcal{H}_{0}^{s}(\Omega)\rightarrow L^{2}(\Omega)$
we can consider its adjoint operator with respect to the scalar product
$\mathcal{B}^{a}$,
\[
i^{*}:L^{2}(\Omega)\rightarrow\mathcal{H}_{0}^{s}(\Omega).
\]
The composition $(i^{*}\circ i)_{a}:\mathcal{H}_{0}^{s}(\Omega)\rightarrow\mathcal{H}_{0}^{s}(\Omega)$
is selfadjoint, compact, injective with dense image in $\mathcal{H}_{0}^{s}(\Omega)$
and it holds
\begin{equation}
\mathcal{B}^{a}\left((i^{*}\circ i)_{a}u,v\right)=\mathcal{E}\left((i^{*}\circ i)_{a}u,v\right)+\int_{\Omega}au(i^{*}\circ i)_{a}v=\int_{\Omega}uv.\label{eq:ii*}
\end{equation}
We call $\varphi^{a}\in\mathcal{H}_{0}^{s}(\Omega)$ an eigenfunction
of $\left(\left(-\Delta\right)^{s}+a\right)$ corresponding to the
eigenvalue $\lambda^{a}$ if 
\[
\mathcal{E}(\varphi^{a},v)+\int_{\mathbb{R}^{n}}a\varphi^{a}vdx=\lambda^{a}\int_{\mathbb{R}^{n}}\varphi vdx\ \ \forall v\in\mathcal{H}_{0}^{s}(\Omega).
\]

Notice that if $\varphi_{k}^{a}\in\mathcal{H}_{0}^{s}(\Omega)$ is
an eigenfunction of the fractional Laplacian with eigenvalue $\lambda_{k}^{a}$,
then $\varphi_{k}^{a}$ is an eigenfunction of $(i^{*}\circ i)_{a}$
with eigenvalue $\mu_{k}^{a}:=1/\lambda_{k}^{a}$. In fact, it holds,
for all $v\in\mathcal{H}_{0}^{s}(\Omega)$
\[
\mathcal{B}^{a}(\varphi_{k}^{a},v)=\lambda_{k}^{a}\int_{\mathbb{R}^{n}}\varphi_{k}^{a}vdx=\int_{\mathbb{R}^{n}}\lambda_{k}^{a}\varphi_{k}^{a}vdx=\mathcal{B}^{a}\left(\lambda_{k}^{a}(i^{*}\circ i)_{a}\varphi_{k}^{a},v\right),
\]
thus $(i^{*}\circ i)_{a}\varphi_{k}^{a}=1/\lambda_{k}^{a}\varphi_{k}^{a}$. 

We recall, also, that $\left(\left(-\Delta\right)^{s}+a\right)$ admits
an ordered sequence of eigenvalues 
\[
0<\lambda_{1}^{a}<\lambda_{2}^{a}\le\lambda_{3}^{a}\le\dots\le\lambda_{k}^{a}\le\dots\rightarrow+\infty.
\]
and the eigenvalues $\lambda_{k}^{a}$ depend continuously on $a$. 

For $b\in C^{0}(\Omega)$ with $\|b\|_{L^{\infty}}$ small enough
consider $\mathcal{B}^{a+b}$ and $(i^{*}\circ i)_{a+b}$ and set
\begin{equation}
B_{b}:=\mathcal{B}^{a+b}\text{ and }E_{b}:=(i^{*}\circ i)_{a+b}.\label{eq:notation}
\end{equation}

We want to apply the abstract Theorem \ref{thm:astratto-1} to the
operator $E_{b}$ to check when a perturbation $b$ preserve the multiplicity
of an eigenvalue $\mu_{k}^{a}$. Thus, in light of previous consideration,
we get the persistence result for the operator $\left(\left(-\Delta\right)^{s}+a\right)$.
Since we endowed $\mathcal{H}_{0}^{s}(\Omega)$ with the scalar product
$\mathcal{B}^{a}=B_{0}$, to check condition (\ref{eq:spezzamentoastratto-1})
we need to compute $B_{0}\left(E'(0)[b]u,v\right)$.

By the identity (\ref{eq:ii*}), differentiating along the coefficient
$a(x)$ we get (see \cite[Lemma 20]{FGMP}) we have that $B_{0}\left(E'(0)[b]u,v\right)+B'(0)[b]\left(E_{0}u,v\right)=0$,
so, by (\ref{eq:scalar}) and by direct computation (see also \cite[Remarks 21 and 22]{FGMP})
\[
-B_{0}\left(E'(0)[b]u,v\right)=B'(0)[b]\left(E_{0}u,v\right)=\int_{\Omega}b(E_{0}u)v=\int_{\Omega}b\left[(i^{*}\circ i)_{a}u\right]v.
\]
So, if f $\mu^{a}$ is an eigenvalue of the map $E_{0}=(i^{*}\circ i)_{a}$
with multiplicity $\nu>1$, and $\varphi_{1}^{a},\dots,\varphi_{\nu}^{a}$
are its $L^{2}$-orthonormal eigenvectors we get
\[
\left(B'(0)[b]E_{0}\varphi_{i}^{a},\varphi_{j}^{a}\right)=\int_{\Omega}bE_{0}(\varphi_{i}^{a})\varphi_{j}^{a}=-\mu^{a}\int_{\Omega}b\varphi_{i}^{a}\varphi_{j}^{a},
\]
 for all $i,j=1,\dots,\nu$. In this case, considering (\ref{eq:spezzamentoastratto-1}),
we have to deal with 
\[
b\mapsto\gamma_{ij}(b)=\gamma_{ij}:=\int_{\Omega}b(x)\varphi_{i}^{a}\varphi_{j}^{a}d\sigma.
\]

\section{The case of Problem (\ref{eq:Pb-coeff2})}

As in the previous section, we want to see how equation (\ref{eq:spezzamentoastratto-1})
translates in the setting of Problem (\ref{eq:Pb-coeff1}). Since
we assumed $a>0$ on $\bar{\Omega}$, we endow the space $L^{2}(\Omega)$
with scalar product and norm given, respectively, by
\[
\langle u,v\rangle_{L^{2}}=\int_{\Omega}auv;\ \ \ \ \ \|u\|_{L^{2}}^{2}=\int_{\Omega}au^{2},
\]
while on $\mathcal{H}_{0}^{s}$ we consider the usual scalar product
$\mathcal{E}(u,v)$. Again we consider the embedding $i:\mathcal{H}_{0}^{s}\rightarrow L^{2}$
and its adjoint operator $i^{*}:L^{2}\rightarrow\mathcal{H}_{0}^{s}$.
Then we have 
\[
\mathcal{E}((i^{*}\circ i)_{a}v,u)=\int_{\Omega}auv\ \ \forall u,v\in\mathcal{H}_{0}^{s}.
\]
As before, the map $(i^{*}\circ i)_{a}$ is selfadjoint, continuous
and compact from $\mathcal{H}_{0}^{s}$ in itself, and if $\varphi^{a}$
is an eigenfunction with eigenvalue $\lambda^{q}$ for the problem
(\ref{eq:Pb-coeff1}), then it is also an eigenfunction for $(i^{*}\circ i)_{a}$
associated to the eigenvalue $\mu^{a}=1/\lambda^{a}$.

In this case we have to compute $\mathcal{E}(E'(0)[b]u,v)$. This
can be computed directly (see \cite[Lemma 26]{FGMP} and we have 
\[
\mathcal{E}(E'(0)[b]u,v)=\int_{\Omega}buv.
\]
So, also in this case, if f $\mu^{a}$ is an eigenvalue of the map
$(i^{*}\circ i)_{a}$ with multiplicity $\nu>1$, and $\varphi_{1}^{a},\dots,\varphi_{\nu}^{a}$
are its $L^{2}$-orthonormal eigenvectors, in the end we have to consider
the same function
\[
b\mapsto\gamma_{ij}(b)=\gamma_{ij}:=\int_{\Omega}b(x)\varphi_{i}^{a}\varphi_{j}^{a}d\sigma.
\]

In the next section we will check the conditions on $b\mapsto\gamma_{ij}(b)$
to prove Thm \ref{thm:coeff}.

\section{Proof of Theorem \ref{thm:coeff}}

The first part of the Theorem is the translation of Theorem \ref{thm:astratto-1}
in our setting, and in the previous two sections we showed that both
for Problem (\ref{eq:Pb-coeff1}) and for Problem (\ref{eq:Pb-coeff2})
the operator $b\mapsto\gamma_{ij}(b)<T'(0)[b]x_{j}^{0},x_{i}^{0}>_{X}$
is $\gamma_{ij}=\int_{\Omega}b(x)\varphi_{i}^{a}\varphi_{j}^{a}d\sigma$,
so the set $\mathscr{I}$ of the $b$ near $0$ such that an eigenvalue
$\lambda_{0}$ maintains the same multiplicity is a subset of $\mathscr{H}=\left\{ b\in C^{1}(\Omega)\ :\ \gamma_{ij}(b)=\rho\textrm{Id}\text{ for some }\rho\neq0\right\} $
and it is a smooth sub-manifold of $C^{1}(\Omega)$ if the set $\left\{ \textrm{Id},\text{\ensuremath{\left(\gamma_{ij}(b)\right)_{ij}\text{ for }}}b\in C^{1}(\Omega),\|b\|_{C^{1}}\text{ small}\right\} $
generates all the symmetric $\nu\times\nu$ matrices. It remains to
show that this last condition is fulfilled when $\nu=2$.

In particular, suppose that $\lambda_{0}$ is an eigenvalue for Problem
(\ref{eq:Pb-coeff1}) or for Problem (\ref{eq:Pb-coeff2}) with multiplicity
$\nu=2$. Let $\varphi_{1}$ and $\varphi_{2}$ be the two $L^{2}$-orthogonal
eigenfunctions relative to $\lambda_{0}$. We want to show that
\[
b\mapsto\left(\int_{\Omega}b\varphi_{i}\varphi_{j}\right)_{ij=1,2}
\]
generates all the symmetric $2\times2$ matrices. To do so, it is
sufficient to prove that 
\[
b\mapsto\left(\int_{\Omega}b\varphi_{1}^{2},\int_{\Omega}b\varphi_{2}^{2},\int_{\Omega}b\varphi_{1}\varphi_{2}\right)
\]
generates $\mathbb{R}^{3}$. Let us suppose, by contradiction, that
there exists a $v=(v_{1},v_{2},v_{3})\neq0$ which is orthogonal to all
$\gamma(b)$. Thus it holds
\[
0=v_{1}\int_{\Omega}b\varphi_{1}^{2}+v_{2}\int_{\Omega}b\varphi_{2}^{2}+v_{3}\int_{\Omega}b\varphi_{1}\varphi_{2}=\int_{\Omega}b(v_{1}\varphi_{1}^{2}+v_{2}\varphi_{2}^{2}+v_{3}\varphi_{1}\varphi_{2})
\]
for all $b\in C^{1}$. This would imply $v_{1}\varphi_{1}^{2}+v_{2}\varphi_{2}^{2}+v_{3}\varphi_{1}\varphi_{2}=0$
almost everywhere on $\Omega$. At this point showing that $\left\{ \varphi_{i}\varphi_{j}\right\} _{ij=1,2}$
are independent as functions on $\Omega$, ends the proof. We follow
the strategy of Micheletti Lupo \cite{DL}, using as a crucial tool
the following result (\cite[Teorema 1.4]{FF}).
\begin{lem}
\label{lem:uncont}Let $u\in D^{s,2}(\mathbb{R}^{n})$ be a weak solution
to (\ref{eq:Pb-coeff1}) or (\ref{eq:Pb-coeff2}) in a bounded domain
$\Omega$ with $s\in(0,1)$ with \textup{$a(x)$ a $C^{1}$ function.}
If $u\equiv0$ on a set $E\subset\Omega$ of positive measure, then
$u\equiv0$ in $\Omega$.
\end{lem}

Let us set $\tau=\varphi_{1}$ and $t=\varphi_{2}$, to simplify notation.
We know that $\tau$ and $t$ are independent functions on $\Omega$.
So also $\tau^{2}$ and $\tau t$, $\tau t$ are $t^{2}$ are. We
want to rule out that, for some $A,B\in\mathbb{R}$, 
\begin{equation}
\tau^{2}=At^{2}+B\tau t.\label{eq:dip}
\end{equation}
By Lemma \ref{lem:uncont} the an eigenvalue vanish at most on a zero
measure set on $\Omega$. If $x$ is such $\tau(x)\neq0$ we can divide
(\ref{eq:dip}) by $\tau^{2}$ and solve, obtaining
\[
\frac{t}{\tau}=\left\{ \begin{array}{c}
c_{1}:=\frac{-B+\sqrt{B^{2}+4A}}{2}\\
c_{2}:=\frac{-B-\sqrt{B^{2}+4A}}{2}
\end{array}\right..
\]
So, there exists a set $E\subset\Omega$ such that
\[
t=\left\{ \begin{array}{cc}
c_{1}\tau & \text{on }E\\
c_{2}\tau & \text{ \text{on }}\Omega\smallsetminus E
\end{array}\right..
\]
At least one set between $E$ and $\Omega\smallsetminus E$ has positive
measure. So, we can suppose that $E$ has it. At this point we construct
an eigenfunction $\varphi_{2}-c_{1}\varphi_{1}$ which is zero on
$E$. This contradicts Lemma \ref{lem:uncont} , so $\varphi_{1}^{2},\varphi_{2}^{2}$
and $\varphi_{1}\varphi_{2}$ are independent as functions on $\Omega$
and we had completed the proof.

\section{Proof of Theorem \ref{thm:dom}}

As anticipated in the introduction, we consider a perturbed domain
as $\Omega_{\psi}:=(I+\psi)\Omega,$ with $\psi\in C^{1}(\mathbb{R}^{n},\mathbb{R}^{n})$,
$\|\psi\|_{C^{1}}$ small enough to ensure that $(I+\psi)$ is invertible.
We denote $J_{\psi}$ as the Jacobian determinant of the mapping $I+\psi$. 

By the change of variables given by the mapping $(I+\psi)$, and denoted
$\tilde{u}(\xi):=u(\xi+\psi(\xi))$, we obtain the bilinear form $\mathcal{B}_{s}^{\psi}$
on $\mathcal{H}_{0}^{s}(\Omega)$ defined in the following formula:

\begin{multline}
\mathcal{E}_{s}^{\Omega_{\psi}}(u,v)=\frac{1}{2}\int_{\mathbb{R}^{n}}\int_{\mathbb{R}^{n}}\frac{(u(x)-u(y))(v(x)-v(y))}{|x-y|^{n+2s}}dxdy\\
=\frac{1}{2}\int_{\mathbb{R}^{n}}\int_{\mathbb{R}^{n}}\frac{(\tilde{u}(\xi)-\tilde{u}(\eta))(\tilde{v}(\xi)-\tilde{v}(\eta))}{|\xi-\eta+\psi(\xi)-\psi(\eta)|^{n+2s}}J_{\psi}(\xi)J_{\psi}(\eta)d\xi d\eta\\
=:\frac{1}{2}\mathcal{B}_{s}^{\psi}(\tilde{u},\tilde{v}).
\end{multline}
Here $\tilde{u},\tilde{v}\in\mathcal{H}_{0}^{s}(\Omega)$ and $u,v\in\mathcal{H}_{0}^{s}(\Omega_{\psi})$.
Notice that $\mathcal{B}_{s}^{0}(\tilde{u},\tilde{v})=\mathcal{E}_{s}^{\Omega}(\tilde{u},\tilde{v})$.

To simplify notation, we define the map
\begin{align*}
\gamma_{\psi} & :\mathcal{H}_{0}^{s}(\Omega_{\psi})\rightarrow\mathcal{H}_{0}^{s}(\Omega);\\
\gamma_{\psi}(u) & :=\tilde{u}(\xi)=u(\xi+\psi(\xi)).
\end{align*}
We recall that the map $\gamma_{\psi}$ is invertible since $\|\psi\|_{C^{1}}$
is small. 

As before, given a bounded domain $D$, we consider the embedding
$i:\mathcal{H}_{0}^{s}(D)\rightarrow L^{2}(D)$ and its adjoint operator
$i^{*}$ with respect to the scalar product $\mathcal{E}_{s}^{D}$.
Again, the composition $E_{D}:=(i^{*}\circ i)_{D}:\mathcal{H}_{0}^{s}(D)\rightarrow\mathcal{H}_{0}^{s}(D)$
is a selfadjoint compact operator and
\begin{equation}
\mathcal{E}_{s}^{D}\left(E_{D}u,v\right)=\int_{D}uv.\label{eq:ii*-0}
\end{equation}
In addition, if $\varphi_{k}\in\mathcal{H}_{0}^{s}(D)$ is an eigenfunction
of the fractional Laplacian with eigenvalue $\lambda_{k}$, then it
is also an eigenfunction of $E_{D}=(i^{*}\circ i)_{D}$ with eigenvalue
$\mu_{k}:=1/\lambda_{k}$.

Now, on $\Omega_{\psi}$, we consider $E_{\psi}:=E_{\Omega_{\psi}}$
and we recast (\ref{eq:ii*}) as 
\[
\mathcal{B}_{s}^{\psi}(\gamma_{\psi}E_{\psi}u,\tilde{v})=\mathcal{E}_{s}^{\Omega_{\psi}}\left(E_{\psi}u,v\right)=\int_{\Omega_{\psi}}uv=\int_{\Omega}\tilde{u}\tilde{v}J_{\psi},
\]
and, set 
\[
T_{\psi}=\gamma_{\psi}E_{\psi}\gamma_{\psi}^{-1}\tilde{u},
\]
we have, for $\tilde{u},\tilde{v}\in\mathcal{H}_{0}^{s}(\Omega)$
\[
\mathcal{B}_{s}^{\psi}(T_{\psi}\tilde{u},\tilde{v})=\int_{\Omega}\tilde{u}\tilde{v}J_{\psi}.
\]
We want to apply Theorem \ref{thm:astratto-1} to the selfadjoint
compact operator $T_{\psi}:\mathcal{H}_{0}^{s}(\Omega)\rightarrow\mathcal{H}_{0}^{s}(\Omega)$.

One has, by direct computation, that
\begin{equation}
\left(\mathcal{B}_{s}^{\psi}\right)'(0)[\psi](T_{0}\tilde{u},\tilde{v})+\mathcal{B}_{s}^{0}(T_{\psi}'(0)[\psi]\tilde{u},\tilde{v})=\int_{\Omega}\tilde{u}\tilde{v}\mathrm{div}\psi.\label{eq:calBprim}
\end{equation}
At this point, we use the results of \cite[Lemma 15 and Corollary 16]{FGMP}
(see also \cite[Thm 1.3]{DFW}), to obtain that, if $\varphi_{i},\varphi_{j}\in\mathcal{H}_{0}^{s}(\Omega)$
are two eigenfunctions with the same eigenvalue $\lambda_{0}$ for
the fractional laplacian (in other words, such that $T_{0}\varphi_{i}=\frac{1}{\lambda_{0}}\varphi_{i}$,
and $T_{0}\varphi_{j}=\frac{1}{\lambda_{0}}\varphi_{j}$), it holds
\[
\left(\mathcal{B}_{s}^{\psi}\right)'(0)[\psi](T_{0}\varphi_{i},\varphi_{j})=-\frac{\Gamma^{2}(1+s)}{\lambda_{0}}\int_{\partial\Omega}\frac{\varphi_{i}}{\delta^{s}}\frac{\varphi_{j}}{\delta^{s}}\,\psi\cdot N\,d\sigma+\int_{\Omega}\varphi_{i}\varphi_{j}\mathrm{div}(\psi)dx
\]
and, by (\ref{eq:calBprim})
\[
\mathcal{B}_{s}^{0}(T_{\psi}'(0)[\psi]\varphi_{i},\varphi_{j})=\frac{\Gamma^{2}(1+s)}{\lambda_{0}}\int_{\partial\Omega}\frac{\varphi_{i}}{\delta^{s}}\frac{\varphi_{j}}{\delta^{s}}\,\psi\cdot N\,d\sigma,
\]
 so the operator in formula (\ref{eq:spezz-dom}) in Thm \ref{thm:dom}
is indeed
\[
\psi\mapsto\gamma_{ij}(\psi):=\int_{\partial\Omega}\frac{\varphi_{i}}{\delta^{s}}\frac{\varphi_{j}}{\delta^{s}}\psi\cdot N
\]
as claimed, and the proof of theorem follows.
\begin{rem}\label{dom}\rm
We notice that, also for an eigenvalue of multiplicity $\nu=2$, repeating
the same strategy of the proof of Thm \ref{eq:spezz-coeff}, one could
construct an eigenvalue $\bar{\varphi}$ for which $\frac{\bar{\varphi}}{\delta^{s}}=0$
on a subset of the boundary $\partial\Omega$ which has positive measure.
We remark that $\frac{\bar{\varphi}}{\delta^{s}}$ plays the role
of $\partial_{N}\bar{\varphi}$ in the local case, and if $\bar{\varphi}$
is an eigenvalue for the local laplacian with Dirichlet boundary condition
for which $\partial_{N}\bar{\varphi}=0$ on a set of $\partial\Omega$
of positive measure, then $\bar{\varphi}\equiv0$ on $\Omega$, by
an application of the unique continuation principle, that $\bar{\varphi}\equiv0$
on $\Omega$. Unfortunately, the extension of this result to the fractional
case seems very challenging and, as far as we know, is far from
being proved. 
\end{rem}

\end{document}